\documentclass[LETTER, 10pt, conference]{ieeeconf}      

\IEEEoverridecommandlockouts                              
\usepackage{graphics} 
\usepackage{epsfig} 
\usepackage{mathptmx} 
\usepackage{times} 
\usepackage{amsmath} 
\usepackage{amssymb}  

\newtheorem{theorem}{Theorem}

\newtheorem{proposition}{Proposition}
\newtheorem{lemma}{Lemma}

\newtheorem{example}{Example}

\title{\LARGE \bf   Parameters  estimation of a noisy  sinusoidal
signal with time-varying amplitude }

\author{Da-yan Liu, Olivier Gibaru and Wilfrid Perruquetti
\thanks{D.Y. Liu is with Paul Painlev\'{e} (CNRS, UMR 8524),  Universit\'{e} de Lille 1, 59650, Villeneuve d'Ascq, France
        {\tt\small  dayan.liu@inria.fr}}
\thanks{O. Gibaru is with Laboratory of Applied Mathematics and Metrology (L2MA), Arts et Metiers ParisTech,   8 Boulevard Louis XIV, 59046 Lille Cedex, France
        {\tt\small  olivier.gibaru@ensam.eu}}
\thanks{W. Perruquetti is with LAGIS (CNRS, UMR 8146), \'{E}cole Centrale de Lille, BP 48,
Cit\'e Scientifique, 59650 Villeneuve d'Ascq, France
        {\tt\small  wilfrid.perruquetti@inria.fr}}
\thanks{D.Y. Liu, O. Gibaru and W. Perruquetti  are with \'{E}quipe Projet Non-A, INRIA Lille-Nord Europe, Parc Scientifique
de la Haute Borne 40, avenue Halley B\^{a}t.A, Park Plaza, 59650
Villeneuve d'Ascq, France
        } }

\begin{document}

\maketitle
\thispagestyle{empty}
\pagestyle{empty}

\begin{abstract}
In this paper, we give  estimators of the frequency, amplitude and
phase  of a  noisy sinusoidal signal with time-varying amplitude by
using the algebraic parametric techniques introduced by Fliess and
Sira-Ram\'{\i}rez. We apply a similar strategy  to estimate these
parameters by using modulating functions method. The convergence of
the noise error part due to a large class of noises is studied to
show the robustness and the stability of these methods. We also show
that the estimators obtained by modulating functions method are
robust to ``large'' sampling period and to non zero-mean noises.
\end{abstract}

\section{INTRODUCTION}

Recent algebraic parametric estimation techniques for linear systems
\cite{mfhsr,garnier,c7} have been extended to various problems in
signal processing (see, e.g., \cite{mexico,mboup,num,Cdc09,JCAM}).
In \cite{Med08,trapero,Trapero3}, these methods are devoted to
estimate the frequency, amplitude and phase of a  noisy sinusoidal
signal with time-invariant amplitude. Let us emphasize that these
methods, which are algebraic and non-asymptotic, exhibit good
robustness properties with respect to corrupting noises, without the
need of knowing their statistical properties (see \cite{ans,shannon}
for more theoretical details). We have shown in \cite{algo} that the
differentiation estimators proposed by  algebraic parametric
techniques can cope with a large class of noises for which the mean
and covariance are polynomials in time. The robustness properties
have already been confirmed by numerous computer simulations and
several laboratory experiments.  In \cite{Fedele, Med08}, modulating
functions methods are  used to estimate unknown parameters of noisy
sinusoidal signals. These methods have similar advantages than
algebraic parametric techniques especially concerning the robustness
of estimations to corrupting noises. The aim of this paper is to
estimate the frequency, amplitude and phase of a noisy time-varying
amplitude  sinusoidal signal by using the previous two methods. We
also show their stability by studying the convergence of the noise
error part due to a large class of noises.

In Section \ref{section2}, we give some notations and useful
formulae. In Section \ref{section3} and  Section \ref{section4}, we
give parameters' estimators  by using respectively  algebraic
parametric techniques and  modulating functions method. In Section
\ref{section5}, the estimators are given in  discrete case. Then, we
study the influence of sampling period on the associated noise error
part due to a class of noises. In Section \ref{section6}, inspired
by \cite{Fedele} a recursive algorithm for the frequency estimator
is given, then some numerical simulations are given  to show the
efficiency and stability of our estimators.


\section{Notations  preliminaries}\label{section2}

Let us denote by $D_T:=\{T\in \mathbb{R}^*_+: [0,T]
\subset \Omega \}$, and $
 w_{\mu,\kappa}(\tau)=(1-\tau)^{\mu} \tau^{\kappa}$ for any  $\tau \in [0,1]$ with $\mu, \kappa \in ]-1,+\infty[$. By
using the Rodrigues formula (see \cite{R35} p.67), we have
\begin{equation}\label{formule_lei}
\begin{split}
\frac{d^i}{d\tau^i}\left\{w_{\mu,\kappa}(\tau)\right\}=(-1)^i
i!\, w_{\mu-i,\kappa-i}(\tau) P_{i}^{\mu-i,\kappa-i}(\tau),
\end{split}
\end{equation}
where $P_i^{\mu-i,\kappa-i}$, $\min(\kappa,\mu) \geq  i \in \mathbb{N}$, is the $i^{th}$ order Jacobi polynomial defined on
$[0,1]$ (see \cite{R35}): $\forall \, \tau \in [0,1],  $
\begin{equation} \label{jacobi}
P_{i}^{\mu-i,\kappa-i}(\tau)=\sum_{s=0}^{i}(-1)^{i-s}\binom{\mu}{s}\binom{\kappa}{i-s}w_{i-s,i}(\tau).
\end{equation}
Then, we have the following lemma.

\begin{lemma}\label{formule_inverse}
Let $f$ be a $\mathcal{C}^{n+1}(\Omega)$-continuous function ($n \in
\mathbb{N}$) 
and
$\Pi_{k,\mu}^{n}$ be a differential operator defined as follows
\begin{equation}
\Pi_{k,\mu}^{n}=\frac{1}{s^{n+1+\mu}}\cdot\frac{d^{n+k}}{ds^{n+k}}\cdot{s^{n}%
},
\end{equation}
where $s$ is the Laplace variable, $k \in\mathbb{N}$ and $-1 <\mu \in \mathbb{R}$.  Then,  the inverse  Laplace transform of $\Pi_{k,\mu}%
^{n}\hat{f}$ where $\hat{f}$ is the laplace transformation of $f$  is given by
\begin{equation} \label{formule}
\begin{split}
&\mathcal{L}^{-1}\left\{  \Pi_{k,\mu}^{n}\hat{f}(s)\right\}  (T) \\
=& T^{n+1+\mu+k} c_{\mu+n,k} \int_{0}^{1} w^{(n)}_{\mu+n,k+n}%
(\tau)\, f(T\tau)d\tau,
\end{split}
\end{equation}
where $T \in D_T$ and
$c_{\mu+n,\kappa}=\frac{(-1)^{\kappa}}{\Gamma(\mu+n+1)}$.
\end{lemma}

In order to prove this lemma, let us recall that the
 $\alpha$-order ($\alpha\in\mathbb{R}^{+}$) Riemann-Liouville integral
(see \cite{FC}) of a  real function $g$
$(\mathbb{R}\rightarrow\mathbb{R})$ is defined by
\begin{equation}\label{}
J^{\alpha}g(t):=\frac{1}{\Gamma(\alpha)}\int_{0}^{t}(t-\tau)^{\alpha-1}%
g(\tau)d\tau.
\end{equation}
The associated Laplace transform is given by
\begin{equation}\label{}
\mathcal{L}\left\{ J^{\alpha }g(t)\right\} (s)=s^{-\alpha}\hat{g}(s),
\end{equation}
 where
$\hat{g}$ denotes the Laplace transform of
$g$.

\noindent\textbf{Proof.} Let us denote $W_{\mu+n,\kappa+n}(t)=(T-t)^{\mu+n} t^{\kappa+n}$ for any $t \in
[0,T]$. Then, by applying
the Laplace transform to the following Riemann-Liouville integral
and doing some classical operational calculations, we obtain
\begin{equation*}%
\begin{split}
&\mathcal{L}\left\{c_{\mu+n,k+n}\int_{0}^{T} W_{\mu+n,k+n}(\tau)
f^{(n)}(\tau)d\tau \right\}\\
=&s^{-(n+1+\mu)} \mathcal{L} \left\{ (-1)^{n+k} \tau^{n+k} f^{(n)}(\tau)\right\}\\
=&s^{-(n+1+\mu)} \frac{d^{n+k}}{d s^{n+k}} \mathcal{L} \left\{ f^{(n)}(\tau)\right\} \\
=&s^{-(n+1+\mu)} \frac{d^{n+k}}{d s^{n+k}} s^n
\hat{f}(s)=\Pi_{k,\mu}^{n}\hat{f}(s).
\end{split}
\end{equation*}
Then,  by substituting  $\tau$ by $T\tau$ we have
\begin{equation}\label{11}
\begin{split}
   & c_{\mu+n,k+n}\int_{0}^{T}
W_{\mu+n,k+n}(\tau) f^{(n)}(\tau)d\tau\\=&
T^{2n+k+\mu+1}c_{\mu+n,k+n}\int_{0}^{1} w_{\mu+n,k+n}(\tau)
f^{(n)}(T\tau)d\tau.
\end{split}
\end{equation}
By using (\ref{formule_lei}), we obtain
$w^{(i)}_{\mu+n,k+n}(0)=w^{(i)}_{\mu+n,k+n}(1)=0$ for $i=0, \cdots,
n-1$. Finally, this proof can be completed  by applying
$n$ times integration by parts to (\ref{11}). \hfill$\Box$\\

\section{Algebraic parametric techniques}\label{section3}

Let $y= x +  \varpi$ be a noisy  observation on a finite time
interval $\Omega \subset \mathbb{R}^+$ of a real valued signal $x$,
where  $\varpi$ is an additive  corrupting  noise and
\begin{equation}\label{}
\forall t \in \Omega, \ x(t)=(A_0+A_1 t) \, \sin(\omega t + \phi)
\end{equation}
 with $A_0 \in \mathbb{R}^*_+$, $A_1 \in \mathbb{R}^*$, $\omega \in \mathbb{R}^*_+$ and
$-\frac{1}{2}\pi <\phi < \frac{1}{2}\pi$. Observe that  $x$ is a
time-variant varying  sinusoidal signal, which is a solution to the
harmonic oscillator equation
\begin{equation}\label{equation_harmonic00}
\forall\, t \in \Omega, \ \  {x}^{(4)}(t)+ 2\omega
^{2}\ddot{x}(t)+\omega ^{4}x(t)=0.
\end{equation}
Then, we can estimate the parameters $\omega$, $A_0$ and $\phi$ by applying  algebraic parametric techniques to (\ref{equation_harmonic00}).

\begin{proposition}\label{p3}
Let $k \in \mathbb{N}$, $-1 <\mu \in \mathbb{R}$ and $T \in D_{T}$ such that
$A_1\int_{0}^{1} \dot{w}_{\mu+4,k+4}
(\tau)\, \sin(\omega T\tau+\phi)d\tau \leq 0$,
 then the
parameter $\omega$  is estimated from the noisy observation $y$ by
\begin{equation}\label{omega}
\tilde{\omega}=\left(
\frac{-B_y+\sqrt{B_y^2-4A_yC_y}}{2A_y}\right)^{\frac{1}{2}},
\end{equation}
where $A_y=T^4\int_{0}^{1} w_{\mu+4,k+4}%
(\tau)\, y(T\tau)d\tau$, $B_y=2 T^2 \int_{0}^{1} \ddot{w}_{\mu+4,k+4}%
(\tau)\, y(T\tau)d\tau$, $C_y=\int_{0}^{1} w^{(4)}_{\mu+4,k+4}%
(\tau)\, y(T\tau)d\tau $,
$w^{(i)}_{\mu+4,k+4}$ is given by (\ref{formule_lei}) with  $i=1,2,4$.
\end{proposition}

\noindent\textbf{Proof.} By applying the Laplace transform to
(\ref{equation_harmonic00}), we get
\begin{equation}
\begin{split}
&s^4 \hat{x}(s)+2\omega^2 s^2\hat{x}(s)+ \omega^4\hat{x}(s)\\=&s^3 x_0
+s^2 \dot{x}_0 +(2\omega^2 {x}_0+\ddot{x}_0)s+(2\omega^2\dot{x}_0+
x^{(3)}_0). \label{TL11}
\end{split}
\end{equation}
 Let us apply   $k+4 \, (k \in
\mathbb{N})$ times derivations to both sides of (\ref{TL11}) with
respect to $s$. By multiplying  the resulting equation by
$s^{-5-\mu}$ with $-1<\mu \in \mathbb{R}$, we get
\begin{equation}
\Pi_{k,\mu}^4 \hat{x}(s) + 2\omega^2 \Pi_{k+2,\mu+2}^2
\hat{x}(s)+\omega^4 \Pi_{k+4,\mu+4}^0\hat{x}(s)=0. \label{TL22}
\end{equation}
Let us apply the inverse Laplace transform to (\ref{TL22}), then by
using Lemma \ref{formule_inverse}, we obtain
\begin{equation*} 
\begin{split}
&\int_{0}^{1} \left(w^{(4)}_{\mu+4,k+4}%
(\tau)  + 2(\omega T)^2  \ddot{w}_{\mu+4,k+4}%
(\tau)\right) x(T\tau)\ d\tau \\&+ (\omega T)^4 \int_{0}^{1} w_{\mu+4,k+4}%
(\tau)\, x(T\tau)d\tau=0.
\end{split}
\end{equation*}
According to (\ref{formule_lei}), we have $w_{\mu+4,k+4}^{(i)}(0)=w_{\mu+4,k+4}^{(i)}(1)$ for
$i=0,\dots,3$. Then by applying integration by parts, we get
\begin{equation*} 
\begin{split}
&\omega^4 \int_{0}^{1} w_{\mu+4,k+4}%
(\tau)\, x(T\tau)+ 2\omega^2  {w}_{\mu+4,k+4}%
(\tau) x^{(2)}(T\tau)\ d\tau \\&+ \int_{0}^{1} w_{\mu+4,k+4}%
(\tau) x^{(4)}(T\tau)d\tau =0.
\end{split}
\end{equation*}
Thus, $\omega^2$ is obtained  by
\begin{equation}\label{omega2}
{\omega^2}=
\frac{-\hat{B}_x \pm\sqrt{\hat{B}_x^2-4 \hat{A}_x \hat{C}_x}}{2\hat{A}_x},
\end{equation}
where
  $\hat{A}_x= \int_{0}^{1} w_{\mu+4,k+4}%
(\tau)\, x(T\tau)d\tau$, $
\hat{B}_x=2 \int_{0}^{1}{w}_{\mu+4,k+4}%
(\tau) x^{(2)}(T\tau)\ d\tau,$
 $\hat{C}_x=\int_{0}^{1} w_{\mu+4,k+4}%
(\tau) x^{(4)}(T\tau)d\tau$. Since ${x}^{(4)}(T\tau)+ 2
\omega^{2}x^{(2)}(T\tau)+ \omega^{4}x(T\tau)=0$  for any $\tau \in
[0,1]$, we get
\begin{equation*}
\begin{split}
&\frac{1}{4}\left(\hat{B}_x^2-4\hat{A}_x \hat{C}_x\right)= \\& \left(\int_{0}^{1} {w}_{\mu+4,k+4}%
(\tau) x^{(2)}(T\tau)\ d\tau+\int_{0}^{1} \omega^{2} {w}_{\mu+4,k+4}%
(\tau) x(T\tau)\ d\tau\right)^2.
\end{split}
\end{equation*}
Observe that $x^{(2)}(T\tau)+\omega^2 x(T\tau)=2\omega A_1 \cos(\omega T t+ \phi)$
for any $\tau \in
[0,1]$. If $\omega A_1 \int_{0}^{1} {w}_{\mu+4,k+4}
(\tau)\,  \cos(\omega T t+ \phi)d\tau= -\frac{A_1}{T}\int_{0}^{1} \dot{w}_{\mu+4,k+4}
(\tau)\, \sin(\omega T\tau+\phi)d\tau \geq 0$, then we get
\begin{equation}
\frac{-\hat{B}_x +\sqrt{\hat{B}_x^2-4 \hat{A}_x \hat{C}_x}}{2\hat{A}_x}={\omega^2}.
\end{equation}
Finally, this proof can be completed by applying integration by parts and  substituting $x$ by $y$ in
the last equation.
\hfill$\Box$

By observing  that $x_0=x(0)=A_0 \sin \phi$, $\dot{x}_0=\dot{x}(0)=A_0
\omega \cos \phi+A_1 \sin \phi$ and $x^{(3)}_0=x^{(3)}(0)=
-\omega^2 \dot{x}_0- 2\omega^2 A_1 \sin \phi$, then we can obtain
$A_0 \cos
\phi=\frac{1}{2\omega^3}\left(x^{(3)}_0+3\omega^2\dot{x}_0\right)$.
Hence, if $-\frac{\pi}{2}< \phi <\frac{\pi}{2}$, then we have
\begin{equation} \label{equations2}
\begin{split}
A_0&=\left(  x_{0}^{2}+  \frac{\left(x^{(3)}_0+3\omega^2\dot{x}_0\right)^2}{4\omega^6}\right)^{\frac{1}{2}%
},\\
 \phi&=\arctan\left( \frac{2\omega^3 x_{0}}{x^{(3)}_0+3\omega^2\dot{x}_0} \right).
\end{split}
\end{equation}
Thus, we need to estimate ${x}_{0}$, $\dot{x}_{0}$ and $x^{(3)}_{0}$
so as to obtain the estimations of
 $A_0$ and $\phi$.

\begin{proposition}\label{p4} Let  $-1 <\mu \in \mathbb{R}$ and $T \in D_{T}$, then the parameters
$A_0$ and $\phi$ are estimated from the noisy observation $y$ and
the estimated value of  $\omega$ given in (\ref{omega}):
\begin{equation}
\begin{split}
\tilde{A}_0&=\left(  \tilde{x}_{0}^{2}+  \frac{\left({{\tilde{x}}}^{(3)}_0+3\tilde{\omega}^2{\tilde{\dot{x}}}_0\right)^2}{4\tilde{\omega}^6}\right)^{\frac{1}{2}%
},\\
 \tilde{\phi}&=\arctan\left( \frac{2\tilde{\omega}^3 x_{0}}{{\tilde{x}}^{(3)}_0+3\tilde{\omega}^2\tilde{\dot{x}}_0} \right),
\end{split}
\end{equation}
where
\begin{equation*}
\begin{split}
\tilde{x}_0 =&\int_0^{1} P_2^{\tilde{\omega}}(\tau)\, y(T\tau) \,d\tau, \quad \tilde{\dot{x}}_0= \frac{1}{T} \int_0^{1} P_3^{\tilde{\omega}} (\tau) \, y(T\tau) \,d\tau, \\ \tilde{{x}}^{(3)}_0 =& \frac{1}{T^3} \int_0^{1} P_4^{\tilde{\omega}} (\tau) \, y(T\tau) \,d\tau -2\tilde{\omega}^2 \tilde{\dot{x}}_0,\\
\end{split}
\end{equation*}
\begin{equation*}
\begin{split}
&\frac{6}{\Gamma(\mu+5)} P_2^{\tilde{\omega}}(\tau)=\sum_{i=0}^3 \binom{3}{i}
\frac{4!\,c_{\mu+i,3-i} }{(4-i)!}w_{\mu+i,3-i}(\tau) \\&+4(\tilde{\omega}
T)^2\sum_{i=0}^2 \binom{3}{i} \frac{ c_{\mu+i+2,3-i}}{(2-i)!}
w_{\mu+i+2,3-i}(\tau)\\&+(\tilde{\omega} T)^4 c_{\mu+4,3}
w_{\mu+4,3}(\tau),
\end{split}
\end{equation*}
\begin{equation*}
\begin{split}
&\frac{-2}{\Gamma(\mu+6)}P_3^{\tilde{\omega}}(\tau)=
c_{\mu,3} w_{\mu,3}(\tau)+11 c_{\mu+1,2} w_{\mu+1,2}(\tau)\\&+28
c_{\mu+2,1}
w_{\mu+2,1}(\tau)+12 c_{\mu+3,0} w_{\mu+3,0}(\tau)\\
&+2(\tilde{\omega} T)^2\left( c_{\mu+2,3} w_{\mu+2,3}(\tau)+5
c_{\mu+3,2} w_{\mu+3,2}(\tau)\right)\\&+4(\tilde{\omega} T)^2\left( c_{\mu+4,1} w_{\mu+4,1}(\tau)-
c_{\mu+5,0} w_{\mu+5,0}(\tau)\right)\\
&+(\tilde{\omega} T)^4 \left(c_{\mu+4,3}
w_{\mu+4,3}(\tau) -c_{\mu+5,2} w_{\mu+5,2}(\tau)\right),
\end{split}
\end{equation*}
\begin{equation*}
\begin{split}
&\frac{-6}{\Gamma(\mu+8)} P_4^{\tilde{\omega}}(\tau)= \sum_{i=0}^3 \binom{3}{i}
\frac{3!\,c_{\mu+i,3-i} }{(3-i)!}w_{\mu+i,3-i}(\tau) \\&+2(\tilde{\omega}
T)^2\sum_{i=0}^1 \binom{3}{i} c_{\mu+i+2,3-i}
w_{\mu+i+2,3-i}(\tau)\\
&+(\tilde{\omega} T)^4 \sum_{i=0}^3
\binom{3}{i}(-1)^i i! c_{\mu+4+i,3-i}w_{\mu+4+i,3-i}(\tau).
\end{split}
\end{equation*}
\end{proposition}

\noindent\textbf{Proof.} In order to estimate  $x_0$, we apply the
following operator
$\Pi_1=\displaystyle\frac{1}{s^{\mu+5}}\cdot\frac{d^3}{ds^3}$ to
(\ref{TL11}) with $-1<\mu \in \mathbb{R}$, which  annihilates each
terms containing $x^{(i)}_0$ for $i=1,2,3$. Then, by using the
Leibniz formula, we get
\begin{equation*}
\begin{split}
\frac{6}{s^{\mu+5}} x_0& = \sum_{i=0}^3 \binom{3}{i} \frac{4!}{(4-i)!}\frac{1}{s^{\mu+1+i}}
\hat{x}^{(3-i)}(s)\\ +2&\omega^2\sum_{i=0}^2 \binom{3}{i}
\frac{2!}{(2-i)!}\frac{1}{s^{\mu+3+i}}
\hat{x}^{(3-i)}(s)+\frac{\omega^4}{s^{\mu+5}}\hat{x}^{(3)}(s).
\end{split}
\end{equation*}
Let us express the last equation in the time domain. By
 denoting  $T$ as the length of the estimation time window we have
 \begin{equation*}
 \begin{split}
\frac{6 T^{\mu+4}}{\Gamma(\mu+5)} x_0&=\int_0^T \sum_{i=0}^3
\binom{3}{i} \frac{4!\,c_{\mu+i,3-i} }{(4-i)!}W_{\mu+i,3-i}(\tau) d\tau\\
+&4\omega^2 \int_0^T\sum_{i=0}^2 \binom{3}{i} \frac{ c_{\mu+i+2,3-i}}{(2-i)!}
W_{\mu+i+2,3-i}(\tau)d\tau\\+&\omega^4 \int_0^T c_{\mu+4,3} W_{\mu+4,3}(\tau)
d\tau.
\end{split}
 \end{equation*}
Hence, by substituting  $\tau$ by $T\tau$, $x$ by  $y$ and taking
the estimation of $\omega$ given in Proposition \ref{p3} we obtain
an estimate for $x_0$. Similarly, we   apply the
 operator
$\Pi_2=\displaystyle\frac{1}{s^{\mu+4}}\cdot\frac{d}{ds}\cdot\frac{1}{s}\cdot\frac{d^2}{ds^2}$
(resp.
$\Pi_3=\displaystyle\frac{1}{s^{\mu+4}}\cdot\frac{d^3}{ds^3}\cdot\frac{1}{s}$)
to (\ref{TL11})  to compute an estimate for $\dot{x}_0$ (resp.
$x^{(3)}_0$). Finally, we get  estimations for $A_0$ and $\phi$ from
relations (\ref{equations2}) by using  the estimations of $x_0$,
${\dot{x}}_0$, $x^{(3)}_0$ and $\omega$. \hfill$\Box$


\section{Modulating functions method}\label{section4}
\begin{proposition}\label{p7}
Let $f$ be a  function belonging to $\mathcal{C}^{4}([0,1])$ which
satisfies the following conditions $f^{(i)}(0)=f^{(i)}(1)$ for
$i=0,\dots,3$. Assume that
$A_1\int_{0}^{1} \dot{f}
(\tau)\, \sin(\omega T\tau+\phi)d\tau \leq 0$ with $T \in D_{T}$,
 then the
parameter $\omega$  is estimated from the noisy observation $y$ by
\begin{equation}
\tilde{\omega}=\left(
\frac{-B_y+\sqrt{B_y^2-4A_yC_y}}{2A_y}\right)^{\frac{1}{2}},
\end{equation}
where $A_y=T^4\int_{0}^{1} f%
(\tau)\, y(T\tau)d\tau$, $B_y=2 T^2 \int_{0}^{1} \ddot{f}
(\tau)\, y(T\tau)d\tau$, $C_y=\int_{0}^{1} f^{(4)}
(\tau)\, y(T\tau)d\tau$.
\end{proposition}

\noindent\textbf{Proof.} Recall that ${x}^{(4)}(T\tau)+ 2
\omega^{2}x^{(2)}(T\tau)+ \omega^{4}x(T\tau)=0$ for any $\tau \in
[0,1]$. As $f$ is continuous on $[0,1]$, then we have
\begin{equation*}
\begin{split}
&\int_0^1 f(\tau)  {x}^{(4)}(T\tau) d\tau+ 2 \omega^{2} \int_0^1
f(\tau) x^{(2)}(T\tau) d\tau \\ &+\omega^{4} \int_0^1 f(\tau) x(T\tau)
d\tau=0.
\end{split}
\end{equation*}
Then, this proof can be completed similarly to the one of Proposition \ref{p3}.
\hfill$\Box$

\begin{proposition}\label{p8}
Let $f_{i}$ for $i=1,\dots,4$ be four continuous functions defined
on $[0,1]$. Assume that there exists $T \in D_T$ such that the
determinant of the matrix $M_{\omega}=(M^{\omega}_{i,j})_{1 \leq i,j \leq 4}$ is different to zero, where
for $i=1,\dots,4$
\begin{equation*}
\begin{split}
 M_{i,1}^{\omega}&=\int_0^1 f_i(\tau) \sin(\omega T\tau)\,d\tau,  
 M_{i,3}^{\omega}= \int_0^1 f_i(\tau) T\tau\, \sin(\omega T\tau)\,d\tau,   \\
M_{i,2}^{\omega}&=\int_0^1 f_i(\tau) \cos(\omega T\tau)\,d\tau, 
M_{i,4}^{\omega}= \int_0^1 f_i(\tau) T\tau\, \cos(\omega
T\tau)\,d\tau.
\end{split}
\end{equation*}
Then, for any $\phi \in]-\frac{\pi}{2},\frac{\pi}{2}[$ the
estimations of $A_0$, $A_1$ and $\phi$ are given by
\begin{equation} \label{formulae}
\begin{split}
\tilde{A}_i & =\left( \left(\tilde{A_i \cos \phi}\right)^2+\left(\tilde{A_i \sin \phi}\right)^2 \right)^{1/2},\\
\tilde{\phi}  &  =\arctan\left( \frac{\tilde{A_0 \sin
\phi}}{\tilde{A_0 \cos \phi}}\right),
\end{split}
\end{equation}
where the estimates of $A_i \cos \phi$ and $A_i \sin \phi$ for $i=0,1$ are
obtained by solving the following linear system
\begin{equation}
M_{\tilde{\omega}} \left(
\begin{array}{cccc}
  \tilde{A_0 \cos \phi}  \\ \tilde{A_0 \sin \phi} \\ \tilde{A_1 \cos \phi} \\ \tilde{A_1 \sin \phi}
\end{array}
 \right)
= \left(
\begin{array}{cccc}
  I_{f_1}^y  \\ I_{f_2}^y \\ I_{f_3}^y \\ I_{f_4}^y
\end{array}
 \right),
\end{equation}
where $I_{f_i}^y=\int_0^1 f_i(\tau)\, y(T\tau) d\tau$ for
$i=1,\dots,4$, and $\tilde{\omega}$ is the estimate of $\omega$
given by Proposition \ref{p7}.
\end{proposition}

\noindent\textbf{Proof.} Let us take an expansion of $x$
\begin{equation*}
\begin{split}
x(T\tau)=& A_0\cos \phi\, \sin(\omega T\tau)+ A_0 \sin \phi\,
\cos(\omega T\tau) \\ &+ A_1 \cos \phi \, T\tau \sin(\omega T\tau)+ A_1
\sin \phi\, T\tau \cos(\omega T\tau),
\end{split}
\end{equation*}
where $\tau \in [0,1]$, $T \in D_T$.  By multiplying both sides of
the last equation by the continuous functions $f_i$ for $i=1,\dots, 4$
and by integrating the resulting equations  between $0$ and $1$, we obtain
\begin{equation*}\label{}
I_{f_i}^x=A_0\cos \phi M_{i,1}^{\omega} + A_0\sin \phi
M_{i,2}^{\omega} + A_1 \cos \phi M_{i,3}^{\omega} + A_1 \sin \phi
M_{i,4}^{\omega}.
\end{equation*}
Then, it yields the following linear system
\begin{equation*}
M_{\omega} \left(
\begin{array}{cccc}
  A_0 \cos \phi  \\ A_0 \sin \phi \\ A_1 \cos \phi \\ A_1 \sin \phi
\end{array}
 \right)
= \left(
\begin{array}{cccc}
  I_{f_1}^x  \\ I_{f_2}^x \\ I_{f_3}^x \\ I_{f_4}^x
\end{array}
 \right).
\end{equation*}

Since $\det(M_{\omega}) \neq 0$,
we obtain  $A_i \cos \phi$ and $A_i \sin \phi$ for $i=0,1$. Finally, the
proof can be completed by substituting $x$ by $y$ in the so obtained
formulae of $A_i \cos \phi$ and $A_i \sin \phi$. \hfill$\Box$


From now on, we  choose functions  $w_{\mu+n,\kappa+n}^{(n)}$ with
$n \in\mathbb{N}$, $\mu, \kappa \in ]-1,+\infty[$ for the previous
modulating functions. Consequently, the estimate for $\omega$ given
in Proposition \ref{p7}  generalizes the estimate  given in
 Proposition \ref{p3}.


\section{Analysis of the errors due to the noise and the sampling period} \label{section5}

\subsection{Two different sources of errors}
Let us assume now  that $y(t_i)= x(t_i) + \varpi(t_i)$ $(t_i \in
\Omega)$ is a noisy measurement of $x$ in discrete case with an
equidistant sampling period  $T_s$. Since $y$ is a discrete
measurement, we  apply the trapezoidal numerical integration method
to approximate the integrals used  in the previous estimators. Let
$\tau_i= \frac{i}{m}$ and $a_i>0$ for $i=0,\dots, m$ with
$m=\frac{T}{T_s} \in \mathbb{N}^*$ (except for $a_0 \geq 0$ and $a_m
\geq 0$)  be respectively the abscissas and the weights for a given
numerical integration method. Weight $a_0$ (resp. $a_m$) is set to
zero in order to avoid the infinite value at $\tau=0$ when
$-1<\kappa<0$ (resp. $\tau=1$ when $-1<\mu<0$). Let us denote by $q$
 the functions obtained in the integrals of our estimators. Then,
we denote by $I_q^y:=\int_0^1 q(\tau)\, y(T\tau) d\tau.$ Hence,
$I_q^y$ is approximated by $I_q^{y,m}:= \displaystyle\sum_{i=0}^m
\frac{a_i}{m} \, q(\tau_i)\, y(T\tau_i).$ By writing  $y(t_i)=
x(t_i) + \varpi(t_i)$, we get $I_q^{y,m}=I_q^{x,m}+e_q^{\varpi,m},$
where $e_q^{\varpi,m}= \displaystyle\sum_{i=0}^m \frac{a_i}{m} \,
q(\tau_i)\, \varpi(T\tau_i).$ Thus the integral $I_q^y$ is corrupted
by two sources of errors:
\begin{itemize}
  \item  the  numerical error which comes from the  numerical integration method,
  \item the noise error
contributions $e_q^{\varpi,m}$.
\end{itemize}
In the next subsection, we study the choice for the sampling period
so as to reduce the noise error contributions.

\subsection{Analysis of the noise error for different stochastic processes}
We assume in this section that the additive  corruption  noise
$\{\varpi(t_i), t_i \in \Omega\}$ is  a  continuous
 stochastic process satisfying  the following conditions
\begin{description}
  \item[$(C_1):$] for any $s,t \geq 0$, $s\neq t$, $\varpi(s)$ and $\varpi(t)$ are
  independent;
    \item[$(C_2):$] the mean value function of $\{\varpi(\tau), \tau \geq 0\}$ belongs to $\mathcal{L}(\Omega)$;
      \item[$(C_3):$] the  variance function of $\{\varpi(\tau), \tau \geq 0\}$ is bounded on $\Omega$.
\end{description}
Note that  white Gaussian noise and Poisson noise satisfy these conditions.
When the value of $T$ is set, then $T_s \rightarrow 0$ is equivalent to   $m \rightarrow +\infty$.   We are going to show the convergence of the noise error contributions when  $T_s \rightarrow 0$.

\begin{lemma} \label{theorem}
Let $\varpi(t_i)$ be a sequence of $\{\varpi(\tau), \tau \geq 0\}$
with an equidistant sampling period $T_s$, where $\{\varpi(\tau), \tau \geq 0\}$ be a continuous  stochastic process satisfying conditions
$(C_1)-(C_3)$.
Assume that  $q \in \mathcal{L}^2([0,1])$, then we have
\begin{equation}\label{limite1}
\begin{split}
\lim_{m \rightarrow +\infty} E\left[e_q^{\varpi,m}\right]&= \int_0^1 q(\tau) E\left[\varpi(T\tau)\right] d\tau, \\
\lim_{m \rightarrow +\infty} Var\left[e_q^{\varpi,m}\right]&=0.
\end{split}
\end{equation}
\end{lemma}

\noindent\textbf{Proof.} Since $\varpi(t_i)$ is a sequence of
independent random variables $(C_1)$, then by using the properties
of mean value and variance functions  we have
\begin{equation} \label{variance}
\begin{split}
E\left[e_q^{\varpi,m}\right]=&
\frac{1}{m}\sum_{i=0}^m a_i \, q(\tau_i)\,
E\left[\varpi(T\tau_i)\right],\\
Var\left[e_q^{\varpi,m}\right]=&
\frac{1}{m^2}\sum_{i=0}^m a_i^2 \, q^2(\tau_i)\,
Var\left[\varpi(T\tau_i)\right].
\end{split}
\end{equation}
According to $(C_3)$, the variance function of $\varpi$ is bounded.
Then we have
\begin{equation} \label{inegalite}
0\leq\frac{1}{m^2}\sum_{i=0}^m a_i^2 \, q^2(\tau_i)\,
\left|Var\left[\varpi(T\tau_i)\right]\right| \leq
U\frac{a(m)}{m}\sum_{i=0}^m \frac{a_i}{m} \,  q^2(\tau_i),
\end{equation}
where $a(m)=\displaystyle\max_{0\leq i \leq m} a_i$ and
$U=\displaystyle\sup_{0 \leq \tau \leq
1}\left|Var\left[\varpi(T\tau)\right]\right| < +\infty$. Moreover,
since $q \in \mathcal{L}^2([0,1])$ and the mean value function of
$\varpi$ is integrable $(C_2)$, then we have
\begin{equation}
\begin{split} \label{limite}
\lim_{m \rightarrow +\infty} E\left[e_q^{\varpi,m}\right]&= \int_0^1 q(\tau) E\left[\varpi(T\tau)\right] d\tau, \\
\lim_{m \rightarrow +\infty}\sum_{i=0}^m \frac{a_i}{m} \, q^2(\tau_i)&=\int_0^1 q^2(\tau) \,d \tau  <
+\infty.
\end{split}
\end{equation}
As all $a_i$ are bounded,  we have
$U\frac{a(m)}{m}\displaystyle\sum_{i=0}^m \frac{a_i}{m} \,
q^2(\tau_i)= 0.$ This proof  is completed. \hfill $\Box$



\begin{theorem}
With the same conditions given in Lemma \ref{theorem}, we have the following convergence
\begin{equation}\label{}
  e_q^{\varpi,m}  \stackrel{\mathcal{L}^2([0,1])}{\longrightarrow}\int_0^1 q(\tau)\,
E[\varpi(T\tau)] \, d\tau, \ \ \text{when } T_s \rightarrow 0.
\end{equation}
Moreover, if noise $\varpi$  satisfies  the following condition
\begin{description}
  \item[$(C_4):$] $E[\varpi(\tau)]=\displaystyle\sum_{i=0}^{n-1}
{\nu}_i \,\tau^i$ with $n \in \mathbb{N}$ and  ${\nu}_i \in \mathbb{R}$,
\end{description}
and $q \equiv w_{\mu+n,\kappa+n}^{(n)}$ with $ \mu,\kappa \in
]-\frac{1}{2}, +\infty[$,  then we have
\begin{equation}\label{}
 \lim_{m \rightarrow +\infty} E\left[e_q^{\varpi,m}\right]= 0,
\end{equation}
and
\begin{equation}\label{}
e_q^{\varpi,m}  \stackrel{\mathcal{L}^2([0,1])}{\longrightarrow} 0,  \ \ \text{when } T_s \rightarrow 0.
\end{equation}
\end{theorem}

\noindent\textbf{Proof.}  Recall that
$E\left[\left(Y_{m}-c\right)^2\right]=Var\left[Y_{m}\right]+\left(E
\left [Y_{m}\right]-c\right)^2$ for any sequence of random variables
$Y_m$ with $c \in \mathbb{R}$, then by using Lemma \ref{theorem},  $e_q^{\varpi,m}$ converges in mean square to  $\int_0^1 q(\tau)\,
E[\varpi(T\tau)] \, d\tau$
 when $T_s \rightarrow 0$.
If $E[\varpi(\tau)]=\displaystyle\sum_{i=0}^{n-1} {\nu}_i \,\tau^i$
and  $\mu,\kappa \in ]-\frac{1}{2}, +\infty[$, then by using the
Rodrigues formula given by (\ref{formule_lei})  we  obtain
$w_{\mu+n,\kappa+n}^{(n)} \in \mathcal{L}^2([0,1])$ and
 $\int_0^1 w_{\mu+n,\kappa+n}^{(n)}(\tau)\,
E[\varpi(T\tau)] \, d\tau=0$. Hence, this proof is
completed. \hfill
$\Box$

\section{Numerical implementations} \label{section6}
\begin{figure}[!ht]
\centering \epsfig{file=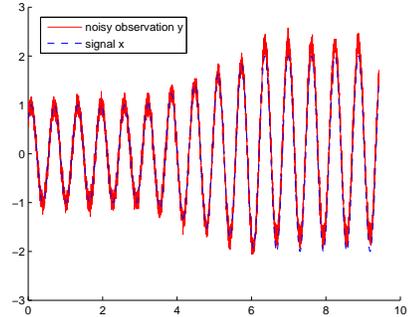,height=4.8cm} \caption{The noisy
observation $y$ and the signal $x$} \label{fig0}
\end{figure}

\begin{figure}[!ht]
\centering \epsfig{file=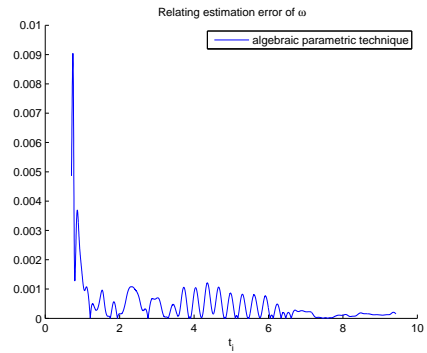,height=4.8cm} \caption{Relating
estimation error of $\omega$} \label{fig1}
\end{figure}

\begin{figure}[!ht]
\centering \epsfig{file=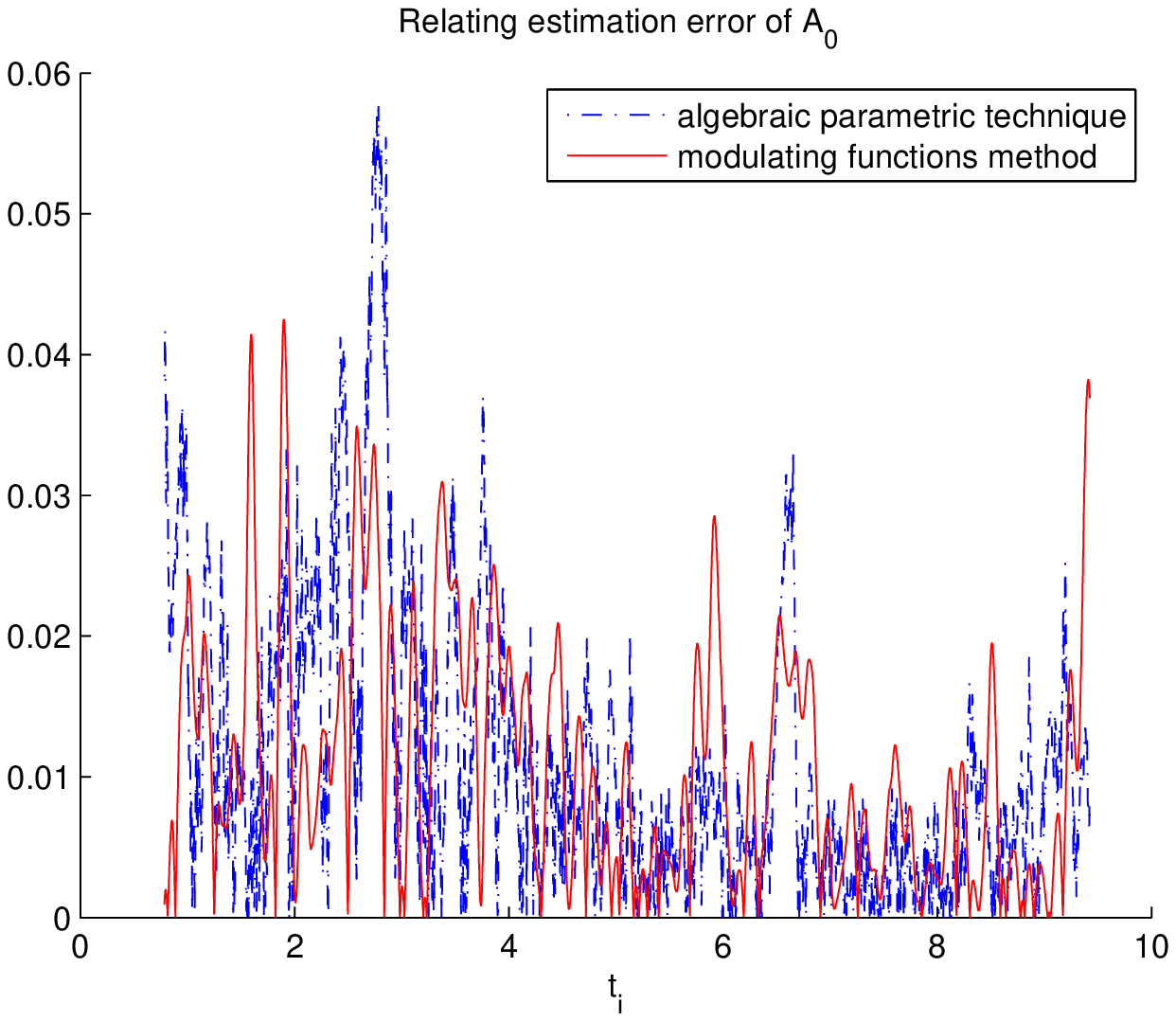,height=4.8cm} \caption{Relating
estimation errors of $A_0$} \label{fig2}
\end{figure}

\begin{figure}[!ht]
\centering \epsfig{file=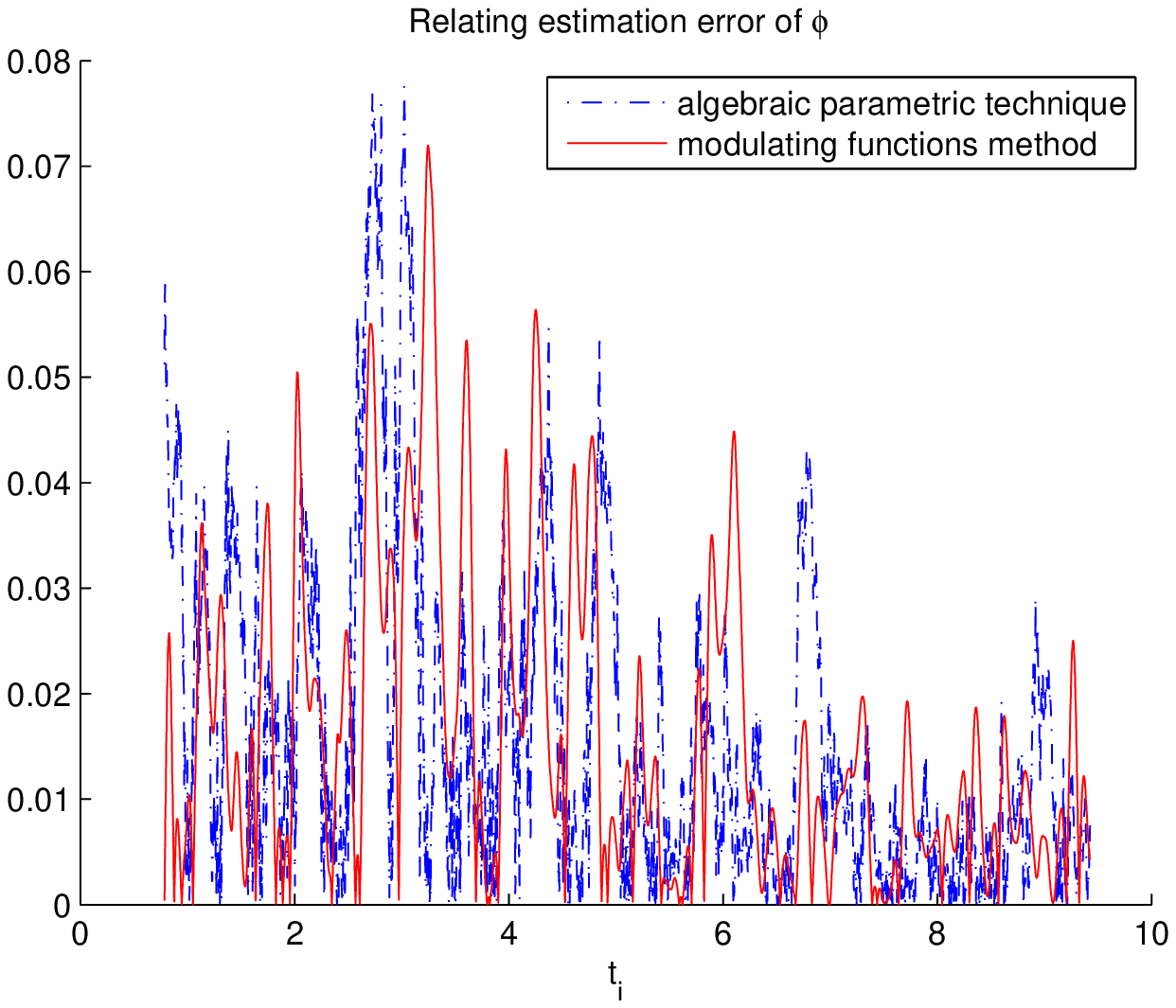,height=4.8cm} \caption{Relating
estimation errors of $\phi$} \label{fig3}
\end{figure}

In our identification procedure, we use a moving integration window.
Hence, the estimate of $\omega$ at $t_i$  is given by Proposition
\ref{p7}  as follows
\begin{equation}\label{discert}
\forall t_i \in \Omega, \ \  \tilde{\omega}^2(t_i)=-
\frac{B_{y_{t_i}}}{2A_{y_{t_i}}}+\frac{\Delta_{y_{t_i}}}{2A_{y_{t_i}}}, \ i=0,1,\dots,
\end{equation}
where $\Delta_{y_{t_i}}=\sqrt{B_{y_{t_i}}^2-4A_{y_{t_i}}C_{y_{t_i}}}$, $A_{y_{t_i}}=T^4  I_f^{y_{t_i},m}$, $B_{y_{t_i}}=2 T^2
I_{\ddot{f}}^{y_{t_i},m}$, $C_{y_{t_i}}=I_{f^{(4)}}^{y_{t_i},m}$ with $y_{t_i} \equiv y(T\cdot+t_i)$.
Note that if $A_{y_{t_i}}=0$, then there is a singular value in $(\ref{discert})$.    If we denote by $\theta_i=\frac{D_{y_{t_i}}}{A_{y_{t_i}}}$ where $D_{y_{t_i}}=-B_{y_{t_i}}$ or $D_{y_{t_i}}=\Delta_{y_{t_i}}$, then
 we can apply the following criterion (see \cite{Fedele})  to improve the estimation of $\omega$
\begin{equation} \label{criterion}
\min_{\theta_i \in \mathbb{R}} J(\theta_i)=\frac{1}{2}\sum_{j=0}^i \nu^{i+1-j}
\left(D_{y_{t_i}}+A_{y_{t_i}} \theta_i\right)^2,
\end{equation}
where $ i=0,1,\dots,$ and  $\nu \in ]0,1]$. The parameter $\nu$   represents a forgetting factor to exponentially discard the ``old'' data in the recursive schema.
The value of $\theta_i$, which minimizes the criterion  $(\ref{criterion})$, is obtained by seeking the value which cancels $\frac{\partial J(\theta_i)}{\partial \theta_i}$. Thus, we get
\begin{equation} \label{theta}
\theta_i=-\frac{\displaystyle\sum_{j=0}^i \nu^{i+1-j}
D_{y_{t_i}}A_{y_{t_i}}}{\displaystyle\sum_{j=0}^i
\nu^{i+1-j}\left(A_{y_{t_i}}\right)^2}.
\end{equation}
Similarly to  \cite{Fedele}, we can get the following recursive algorithm for $(\ref{theta})$
\begin{equation} \label{algorithm}
\theta_{i+1}=\frac{\nu}{\alpha_{i+1}}\left(\alpha_{i}\theta_{i}+D_{y_{t_{i+1}}}A_{y_{t_{i+1}}}\right), \  i=0,1,\dots,
\end{equation}
where $\alpha_{i}=\displaystyle\sum_{j=0}^i
\nu^{i+1-j}\left(A_{y_{t_i}}\right)^2$
Moreover, $\alpha_{i+1}$ can be recursively calculated as follows
$
\alpha_{i+1}=\nu\left(\alpha_{i}+\left(A_{y_{t_i}}\right)^2\right).
$
\begin{example}
According to Section \ref{section5}, we can reduce the noise error
part in our estimations by decreasing the sampling period. Hence,
let $\left(y(t_i)= x(t_i) + c \varpi(t_i)\right)_{i \geq 0}$ be a
generated noise data set  with a small sampling period
$T_s=5\pi\times 10^{-4}$ in the interval $[0,3\pi]$ (see Fig.
\ref{fig0}) where
\begin{equation}\label{signal}
x(t_i)=\left\{
         \begin{array}{rr}
           \sin(10 t_i + \frac{\pi}{4}), & \text{ if } 0 \leq t_i \leq \pi, \\
           \frac{t_i}{\pi}\sin(10 t_i + \frac{\pi}{4}), & \text{ if } \pi < t_i \leq 2 \pi, \\
           2 \sin(10 t_i + \frac{\pi}{4}), & \text{ if } 2 \pi < t_i \leq 3
\pi,
         \end{array}
       \right.
\end{equation}
and noise $c \varpi(x_i)$ is simulated from a zero-mean white
Gaussian $iid$ sequence with $c=0.1$. Hence, the signal-to-noise
ratio $SNR=10\log_{10}\left(\frac{\sum |y(t_i)|^2}{\sum
|c\varpi(t_i)|^2}\right)$ is equal to $SNR=20.8 \text{dB}$.
 In order to estimate the
frequency, by applying  the previous  recursive algorithm we use
Proposition \ref{p3} with $\kappa=\mu=0$, $m=450$ and $\nu=1$. The
relating estimation error is shown in Fig. \ref{fig1}. By using the
estimated  frequency value, we estimate the amplitude and  phase of
the signal   by  applying Proposition \ref{p4} with $\mu=0$, $m=500$
and Proposition \ref{p8} with $m=500$, $f_1 \equiv w_{3,2}$, $f_2
\equiv w_{2,3}$, $f_3 \equiv w_{3,4}$ and $f_4 \equiv w_{4,3}$. The
relating estimation errors  are shown in Fig. \ref{fig2} and Fig.
\ref{fig3}. We can observe that with small  value of $T_s$  the
relating estimation errors are also small.
\end{example}

\begin{figure}[!ht]
\centering \epsfig{file=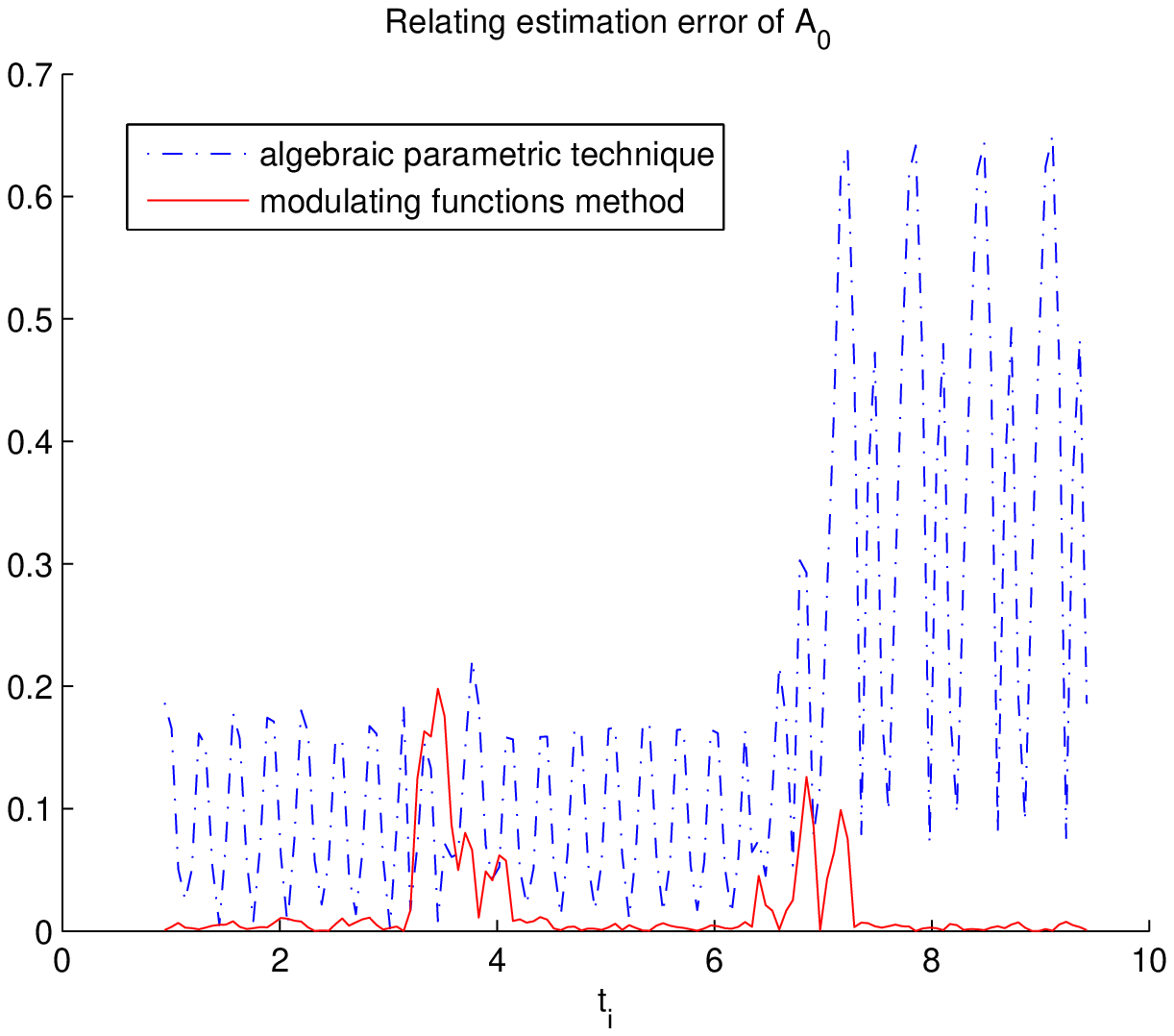,height=4.8cm} \caption{Relating
estimation errors of $A_0$} \label{fig5}
\end{figure}

\begin{figure}[!ht]
\centering \epsfig{file=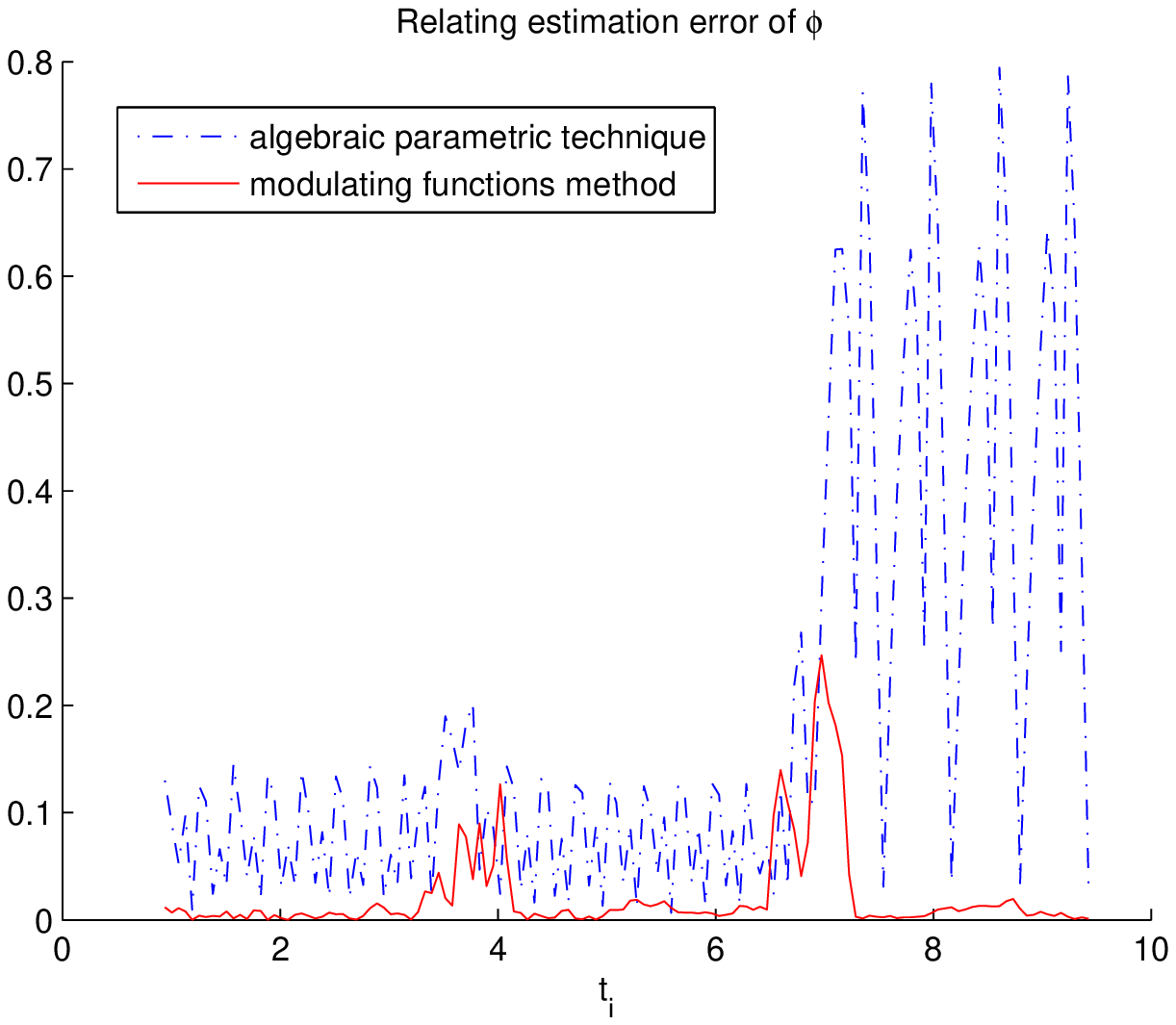,height=4.8cm}
\caption{Relating estimation errors of $\phi$} \label{fig6}
\end{figure}

\begin{example}
In this example, we increase the value of $T_s$ to $T_s=2\pi\times
10^{-2}$ and reduce  the noise level to $c=0.01$. Moreover, we add a
bias term perturbation $\xi=0.25$ in (\ref{signal}) when  $t_i \in
]2 \pi, 3 \pi]$. The estimations of $\omega$ are obtained  by
Proposition \ref{p3} with $\kappa=\mu=0$, $m=12$ and
$\nu=1$. 
The estimations of the amplitude and  phase are given by applying
Proposition \ref{p4} with $\mu=0$, $m=12$ and Proposition \ref{p8}
with $m=15$, $f_1 \equiv w^{(1)}_{3,2}$, $f_2 \equiv w^{(1)}_{2,3}$,
$f_3 \equiv w^{(1)}_{3,4}$ and $f_4 \equiv w^{(1)}_{4,3}$. The
 relating estimation errors  are shown in Fig. \ref{fig5}
and Fig. \ref{fig6}. We can observe that the estimators obtained by
modulating functions method are more robust to the sampling  period
and to the non zero-mean noise than the ones obtained by algebraic
parametric techniques.
\end{example}

\section{CONCLUSIONS AND FUTURE WORKS}
In this paper, two methods are given  to estimate  the frequency,
amplitude and phase  of a  noisy sinusoidal signal with time-varying
amplitude, where the estimates are obtained by using integrals.
There are two types of errors for these estimates: the numerical
error and the noise error part. Then, the convergence in mean square
of the noise error part is studied. A recursive algorithm for
frequency estimator is given. In numerical examples, we show some
comparisons between the two proposed methods. Moreover, these
methods can also be used to estimate the frequencies, the amplitudes
and the phases of two sinusoidal signals from their noisy sum (see
\cite{Trapero3}). The analysis for colored noises will be done in a
future work.

%
%
%

\end{document}